\documentclass[conference]{IEEEtran}
\IEEEoverridecommandlockouts
\usepackage{cite}
\usepackage{amsmath,amssymb,amsfonts,dsfont,balance}
\usepackage{algorithmic}
\usepackage{graphicx,nicefrac}
\usepackage{textcomp}
\usepackage{xcolor}
\usepackage[thinc]{esdiff}
\usepackage{bbm}
\usepackage{stfloats}
\usepackage{balance}

\def\BibTeX{{\rm B\kern-.05em{\sc i\kern-.025em b}\kern-.08em
    T\kern-.1667em\lower.7ex\hbox{E}\kern-.125emX}}

\newtheorem{lemma}{\textbf{Lemma}}
\newtheorem{remark}{\textbf{Remark}}
\newcommand{\longthmtitle}[1]{\mbox{}{\textit{(#1):}}}

\newcommand \bA{\mathbf{A}}
\newcommand \ba{\mathbf{a}}

\newcommand \bI{\mathbf{I}}

\newcommand \br{\mathbf{r}}

\newcommand \bx{\mathbf{x}}
\newcommand \bk{\mathbf{k}}

\newcommand \btau{\boldsymbol{\tau}}
\newcommand \bDelta{\boldsymbol{\Delta}}

\newcommand \bnu{\boldsymbol{\nu}}
\newcommand \bchi{\boldsymbol{\chi}}
\newcommand \bbx{\overline{\mathbf{x}}}
\newcommand \bone{\mathbf{1}}

\newcommand \bzero{\mathbf{0}}

\newcommand \mrss{\mathrm{ss}}

\newcommand \hd{\widehat{d}}
\newcommand \diag{\mathrm{diag}}

\newcommand{\oprocendsymbol}{\hbox{$\bullet$}}
\newcommand{\oprocend}{\relax\ifmmode\else\unskip\hfill\fi\oprocendsymbol}

\DeclareMathOperator*{\argmin}{arg\,min}
\newtheorem{assumption}{Assumption}
 
\begin{document}

\title{Towards Optimal Primary- and Secondary-control Design for Networks with Generators and Inverters\\
{\thanks{This work was supported by the U.S. Department of Energy (DOE) Office of Energy Efficiency and Renewable Energy under Solar Energy Technologies Office (SETO) through the award numbers EE0009025 and 38637 (UNIFI consortium), respectively.}}
}

\author{\IEEEauthorblockN{Manish K. Singh, D. Venkatramanan, and Sairaj Dhople}
\IEEEauthorblockA{\textit{Department of Electrical \& Computer Engineering} \\
\textit{University of Minnesota}\\
Minneapolis, MN USA \\
\{msingh, dvenkat, sdhople\}@umn.edu\vspace{-1em}}
}

\maketitle

\begin{abstract}
For power grids predominantly featuring large synchronous generators (SGs), there exists a significant body of work bridging optimization and control tasks. A generic workflow in such efforts entails: characterizing the steady state of control algorithms and SG dynamics; assessing the optimality of the resulting operating point with respect to an optimal dispatch task; and prescribing control parameters to ensure that (under reasonable ambient perturbations) the considered control nudges the system steady state to optimality. Well studied instances of the aforementioned approach include designing: i) automatic generation control (AGC) participation factors to ensure economic optimality, and ii) governor frequency-droop slopes to ensure power sharing. Recognizing that future power grids will feature a diverse mix of SGs and inverter-based resources (IBRs) with varying control structures, this work examines the different steps of the optimization-control workflow for this context. Considering  a representative model of active power-frequency dynamics of IBRs and SGs, a characterization of steady state is put forth (with and without secondary frequency control). Conditions on active-power droop slopes and AGC participation factors are then derived to ascertain desired power sharing and ensure economically optimal operation under varying power demands.
\end{abstract}

\section{Introduction}
A large and complex interconnected power system can be viewed as  one massive machine with sophisticated controls enabling its secure operation~\cite{Julie_2017}. In most power systems, various parts are owned, operated, and frequently replaced, by loosely coordinating stakeholders. Such flexibility is enabled by an elegant and robust hierarchical control scheme~\cite{COHN1984recollection}. These schemes allow integration of heterogeneous synchronous generators, as long as they comply with certain local dynamic requirements referred to as primary control; and respond to control signals from a supervisory scheme known as secondary control. Being central to the operation of large power grids, these schemes have been actively researched since inception~\cite{Brandt29,cohn1966control,Willems74LFC,Dunlop1975dynamics,Jaleeli_1992,HiskensPaiAGC}. However, there has been a renewed interest in their analysis and design in response to the ongoing transformations in power systems~\cite{Alejandro_2013,JSP_2021,Julija_2015_2,JSP_2022}. These transformations are mainly due to (and to accommodate) the increasing share of inverter-based renewable generation.

Several indicators suggest that inverter-based resources will be commonplace across transmission and distribution networks in future grids~\cite{Ben_2017,NREL_2020}. Vendors of inverter technology are quick to innovate, but incorporating uncountable such devices into the operation and control architecture of power grids is not (and cannot be expected to be) seamless. Utilities and system-operators require easily verifiable benchmarks of anticipated operation and performance. Establishing universal benchmarks is challenged by the different operation timescales, the complex set of nonlinear resource dynamics, and the lack of available models for inverters~\cite{DVR_2022,SVD22Necsys,GrossDorflerReview}. 

We are interested in capturing all salient features of power-system operation without delving deep into the minutiae pertaining to modeling.
\begin{figure}
    \centering
    \includegraphics[width = 1\linewidth]{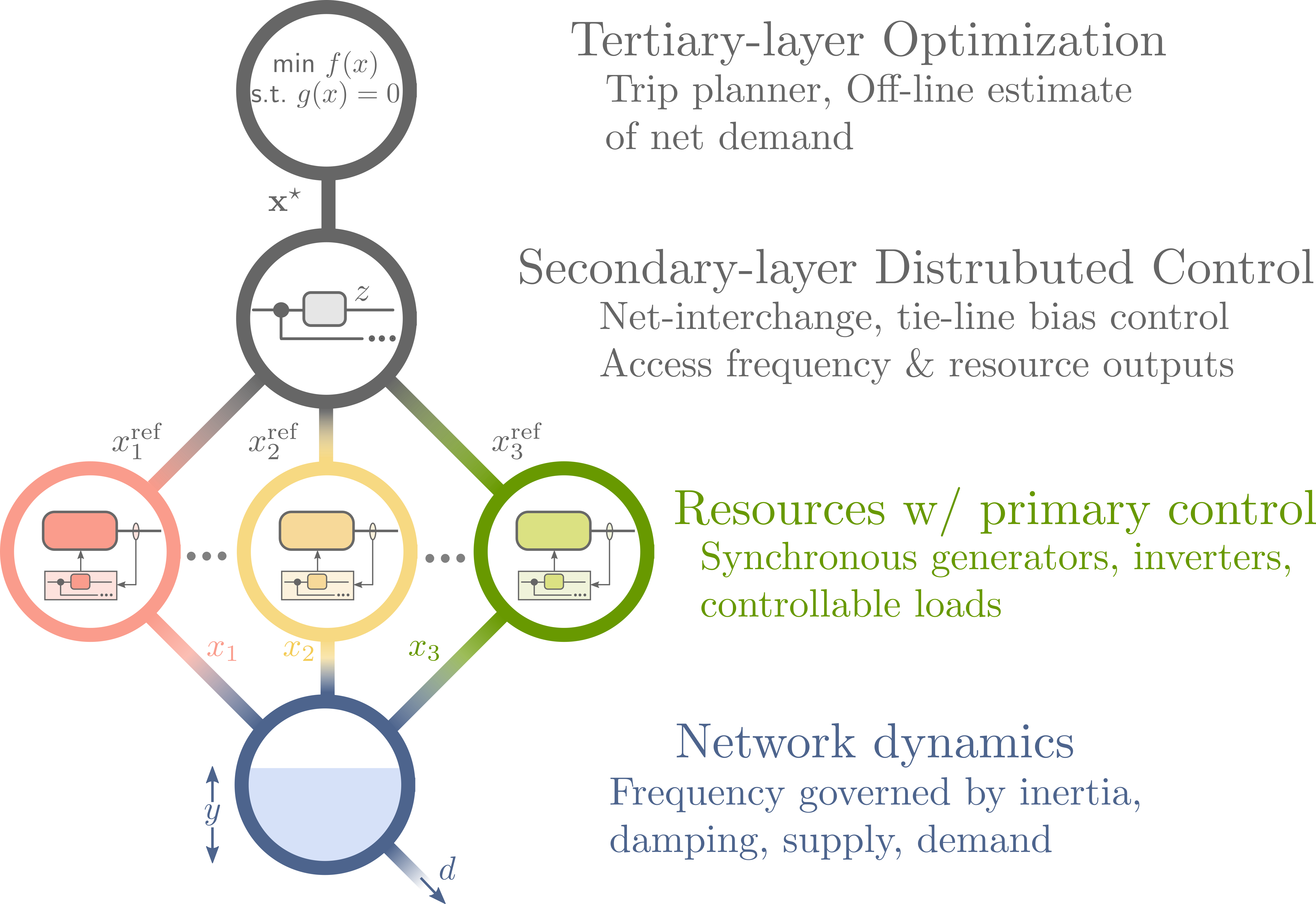}
    \caption{Overview of the considered hierarchical scheme for system operation and control.}
    \vspace{-1em}
    \label{fig:block sch}
\end{figure}
 To that end, we develop a general multi-agent state-space model for a collection of resources (generators and inverters) with outputs (powers) that are intended to supply an external demand. The supply-demand difference is set up to excite the dynamics of a system-wide performance metric (frequency). A perfect match in total supply and demand maps to desired system performance (frequency equals nominal). The resources can respond locally to changes in the performance metric. For higher-level coordination, we have a trip planner that is engaged in optimal resource allocation, as well as a supervisory secondary-control layer that systematically allocates the optimal trajectories to the resources while executing a closed-loop control maneuver to compensate the error in the performance metric. The above architecture mimics all three operational layers of the power system. It is set up with a level of generality to capture a wide range of reported secondary-control schemes as well as macro-level salient features of resource dynamics. See Fig.~\ref{fig:block sch} for an illustration of the overall system architecture. With this generalized state-space model, we obtain the following: 
\begin{itemize}
    \item Steady-state characterization for power system dynamics with and without secondary-control.
    \item Conspicuous dependence of system steady-state on specific design parameters.
    \item Rationale behind conventional design choices for primary- and secondary-control handles.
\end{itemize}

Our development of a generalized system state-space model to mirror fundamentals of power-system architecture is driven by the limited attention bestowed by classical texts on integrated system modeling (and therefore, system-theoretic analysis and design of desirable equilibria and control handles)~\cite{Kundur-1994,Bruce_2010,Glover_2008}. In general, there is an astonishingly vast collection of work on redesigning secondary control, and this has largely mirrored fashionable control synthesis tools du jour~\cite{JSP_2021,Jaleeli_1992,Alejandro_2013,Kothari_2005,JSP_2022_1}. Curiously, there is limited literature on bench-marking the performance of the classical system---as it is understood to be implemented in practice---from the point of view of dynamic and steady-state performance. Notable exceptions to this include~\cite{Baros21AGC} where a connection between dispatch and secondary-control participation factors is established to justify a common design practice for the latter,~\cite{SVD_2020} where the equilibria of the system dynamics are teased out and precisely mapped to power-flow equations, and~\cite{JSP_2022} where a theoretical stability analysis of automatic generation control is presented with a formal treatment of time-scale separation between primary and secondary control.

Of relevance is also a wide body of literature on microgrid control and optimization~\cite{AliMicrogirdTrends,Florian_2016,Davoudi_2020,Davoudi_2018,FlorianHierarchy}. In many such efforts, there is a tendency to mirror the operational layers of the bulk power system. However, given the increased flexibility on offer in synthesizing control and optimization schemes (untethered by adherence to utility and system-operator specifications), several advanced methods involving distributed and decentralized schemes have appeared in the literature. While not directly capturing all features of bulk-system operation (and constraints), these effort, nonetheless, throw up several interesting insights that can be translated.

The remainder of the paper is organized as follows: Section~\ref{sec:modeling} presents the considered abstract dynamic models and overviews their relevant power-system instantiations. Section~\ref{sec:analysis} provides an analysis of the modeled systems, with and without secondary control while particularly emphasizing steady-state characterization. Having elucidated the impact of prominent parameter choices on the system steady state, Section~\ref{sec:design} outlines prudent design choices to harness engineering-centric performance objectives. Summarizing remarks and forward-looking discussions conclude the paper.

\emph{Notation.} Throughout the paper, $\mathds{R}$ denotes the set of real numbers, while $\mathds{R}_{>0}$ and $\mathds{R}_{<0}$ are the subsets of positive and negative real numbers, respectively. Upper- and lower-case boldface letters denote matrices and column vectors. Given a vector $\ba$, its $n$-th entry is denoted by $a_n$; $\diag{(\ba)}$ represents a diagonal matrix with $\ba$ being the principal diagonal. The symbol $(\cdot)^\top$ stands for transposition and $\dot{\bx}$ implies the time-derivative of $\bx$. The symbol $\bone$ denotes a vector of all ones, $\bzero$ denotes a vector or a matrix of all zeros, and $\bI$ is the identity matrix; all with appropriate dimensions.

\section{Modeling}\label{sec:modeling}
For majority of the technical exposition of this paper, we employ a generic notation that considers a resource supply-demand setup with explicit dynamical structure assumed for individual agents, and for a shared supervisory control. The assumed schemes are inspired by prevalent primary- and secondary-control architectures in power systems. The generic structure will allow us to obtain representative results while avoiding cumbersome notational overhangs. Nevertheless, pertinent remarks delineating the mapping from the considered generic model to specific power system instantiations will be suitably provided.

\begin{figure*}
    \centering
    \includegraphics[width = 0.85\linewidth]{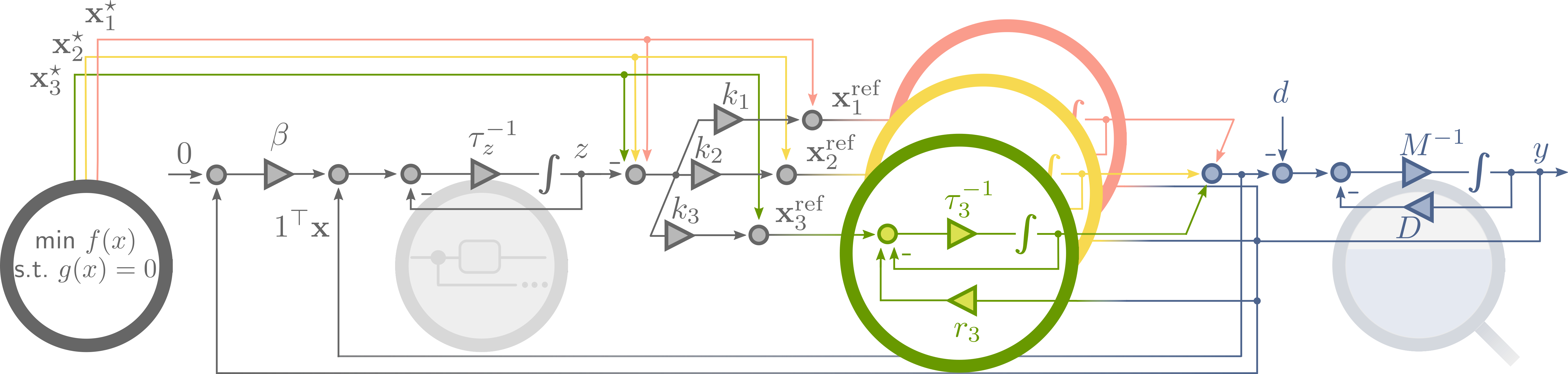}
    \caption{Illustrative high-level block schematic of the hierarchical feedback control system governing the operation of multi-agent power system.}
     \vspace{-1em}
    \label{fig:Block diag}
\end{figure*}

\subsection{System Dynamical Model}
Consider a supply-demand setup of a certain commodity where a demand $d$ is to served collectively by $N$ agents. Denote the supply by vector $\bx\in\mathds{R}^N$ such that $x_n$ is the output of supplier $n$. We will be treating vector $\bx$ as a dynamic state, the governing model of which will shortly follow. Consider a system-wide mismatch indicator $y$ that is a proxy for supply-demand mismatch under steady-state, necessitating specifically $\bone^\top\bx_\mrss=d\iff~y_\mrss=0$, where $\bx_\mrss$ and $y_\mrss$ denote the steady-state values of $\bx$ and $y$, respectively. 

Figure~\ref{fig:block sch} illustrates the overall system we study. It is composed of four subsystems. The first captures the dynamics of the system-wide performance metric, $y$; the second encapsulates the dynamics of the agents (executing control at a primary-control level); the third is a supervisory system effecting control at a secondary level; and the fourth is a long-horizon trip planner implementing a tertiary-level optimization. We discuss the salient dynamics of each subsystem next. 

\subsubsection{System-wide Performance Metric}
The dynamics of state $y$ is assumed to be governed by:
\begin{equation}\label{eq:y_dynamics}
    M \dot{y} +  D y = \bone^\top\bx - d,
\end{equation}
where, $M>0$ and $D>0$ represent cumulative inertia and damping, respectively. 

\subsubsection{Agents \& Primary Control}
The agents are assumed to be first-order integrators with the capacity to respond locally to the mismatch. In particular, they are characterized by dynamics
\begin{align}\label{eq:x_dynamics}
    \diag(\btau)\dot{\bx} &= \bx^{\mathrm{ref}}-\bx + \br y,
\end{align}
where, $\btau \in\mathds{R}^N_{>0}$ captures the time-constants of the agents' response, $\bx_{\mathrm{ref}}\in\mathds{R}^N$ denotes the reference command being provided from the supervisory control layer, and $\br \in\mathds{R}^N_{<0}$ denotes the vector of local regulation gains. Since a positive $y$ is indicative of an oversupply (cf. the steady-state~\eqref{eq:y_dynamics}), local corrective action is engineered for negative feedback by picking entries of $\br$ to be negative.

\subsubsection{Secondary Controller}
The assumed supply dynamics~\eqref{eq:x_dynamics} features local sensing and actuation. Such schemes are known (and will be shown in Section~\ref{sec:analysis}) to be incapable of achieving supply-demand balance, and hence the desired attribute $y_\mrss=0$. Since the demand $d$ is typically an exogenous quantity, a supervisory control can be designed to continually monitor the metric $y$ and supply $\bx$ to generate the time-varying references $\bx^{\mathrm{ref}}$, with an aim of ensuring $y_\mrss=0$ for any constant demand $d$. Denoting by $z$ the state of such a supervisory control, a typical scheme takes the form
\begin{subequations}\label{eq:agc}
\begin{align}
    \tau_z \dot{z} &= -z + \beta y + \bone^\top\bx,\label{seq:zdot} \\
    \bx^\mathrm{ref} &=  \bx^\star + \bk (z - \bone^\top \bx^\star)\label{seq:xref}.
    \end{align}
\end{subequations}
In particular,~\eqref{seq:zdot} represents a generic compensator with time constant, $\tau_{z}$, and secondary-control gain, $\beta<0$. A feed-forward term, $\bone^\top \bx$, is introduced to reduce the burden on the compensator. In~\eqref{seq:xref}, we establish how the references for the agents, $\bx^\mathrm{ref}$ are generated. Particularly, the baseline reference provided by tertiary control, $\bx^\star$, (more on this shortly), is perturbed by a fraction of $z - \bone^\top \bx^\star$ before being relayed to the agents. The idea here being that $z - \bone^\top \bx^\star$ represents the supply needed above and beyond suggested by the off-line trip planner to regulate $y\to0$. The disaggregation of the excess supply is governed by the secondary-control participation factors, $\bk \in \mathds{R}^N$. To ensure that supply adjustments for all agents are constructively aligned towards fulfilling the additional requirement, the participation factors are chosen to be non-negative, i.e., $\bk\geq\bzero$.

\subsubsection{Tertiary Controller}
The tertiary controller generates $\bx^\star$, a constant baseline value of supply that the supervisory control takes as an input. Motivated by economic interests, such a baseline could be set to a minimizer of a resource-allocation problem aimed at serving an anticipated demand $\widehat{d}$. If the cost functions $f_n(\cdot)$ are known for each supplier, the resource-allocation problem can be expressed as
\begin{subequations}
\begin{align}
    \bx^\star=\argmin_{\bchi}~&~\sum_{n=1}^N f_n(\chi_n)\label{eq:RA}\tag{P1}\\
    \mathrm{s.to}~&~\bone^\top\bchi=\widehat d\label{eq:sumx}\\
    &~\bzero\leq\bchi\leq\bbx,\label{eq:limitx}
\end{align}
\end{subequations}
where, \eqref{eq:sumx} ensures supply-demand balance, and~\eqref{eq:limitx} ensures that the supply is bounded within the agents' nominal capacity.

Figure~\ref{fig:Block diag} illustrates a block-schematic representation of the system dynamics~\eqref{eq:y_dynamics}, agent-level local control~\eqref{eq:x_dynamics}, supervisory control~\eqref{eq:agc}, and resource-allocation problem~\eqref{eq:RA}. The power-system parlance for these will be established next.

\subsection{Power-system Instantiation}\label{subsec:PSmodel}
In the context of contemporary power systems, the agents in the system are constituted by synchronous generators and power electronics inverter-based resources (IBRs). The commodity that is collectively supplied by these agents is electrical power. Equation~\eqref{eq:y_dynamics} captures the behaviour of the network frequency arising cumulatively out of either (a) the generator physics or (b) IBR control dynamics~\cite{DVR_2022}. Note that this is traditionally referred to as the swing equation in the context of synchronous machines~\cite{Kundur-1994}, with $\bx$ denoting the actuation input from the mechanical prime mover. Equation~\eqref{eq:x_dynamics} represents the dynamics pertaining to the actuation input, incorporating a droop characteristics in reference command with respect to frequency captured by the gain $\br$~\cite{Bruce_2010}. The secondary layer is also known as automatic generation control (AGC) in classical power system parlance~\cite{Kothari_2005,JSP_2021}. Its dynamics are captured in~\eqref{eq:agc}, where,  the secondary controller output $z$ represents a power command and the quantity $\beta$ is a tunable gain known as the frequency bias (pertaining to a given area). The tertiary control is typically an economic dispatch layer running continuously to optimally utilize the generation assets participating in the power system operation~\cite{SVD_2020}.

\begin{figure*}
    \centering
    \includegraphics[width = 0.85\linewidth]{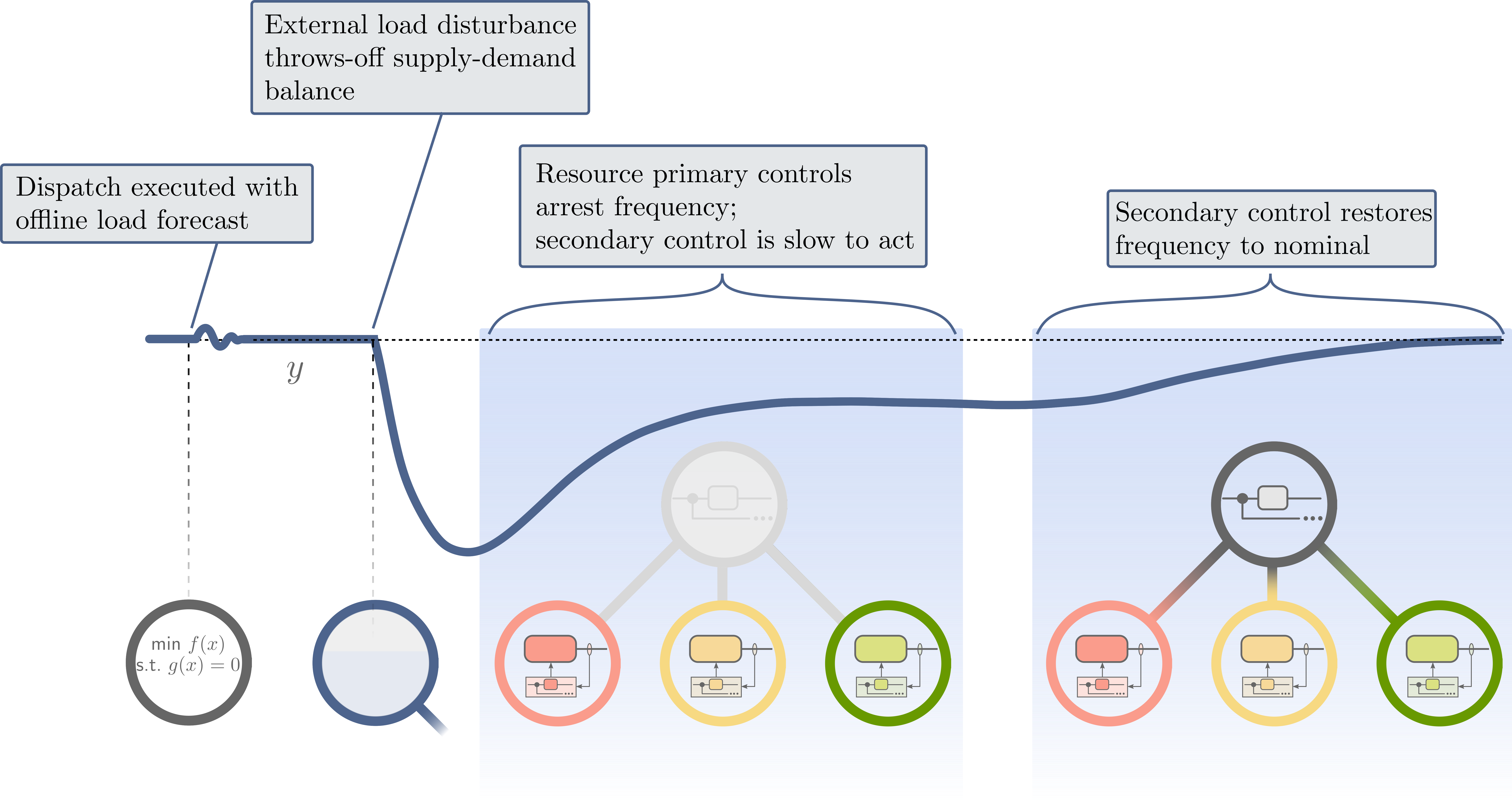}\vspace{-2em}
    \caption{Evolution of system frequency dynamics after disturbance and the associated time-scale separation noted among the hierarchical control layers.}
    \label{fig:Block diag}
     \vspace{-1em}
\end{figure*}

\section{Analysis}\label{sec:analysis}
In this section we analyze the (quasi) steady-state performance of the dynamical system described by~\eqref{eq:y_dynamics}-\eqref{eq:agc}. We first consider the complete system~\eqref{eq:y_dynamics}-\eqref{eq:agc}, and then study a first-order approximation based on time-scale separation.

\subsection{Equilibrium with Secondary Control}
Treating $(\bx^\star,d)$ as exogenous inputs and $(y,\bx,z)$ as states, the dynamics~\eqref{eq:y_dynamics}-\eqref{eq:agc} can be collectively expressed as the linear time-invariant (LTI) system:
\begin{align}\label{eq:statespace}
    \begin{bmatrix}  \dot{y}   \\
         \dot{\bx}  \\
	\dot{z}\end{bmatrix} &=  \begin{bmatrix}  -\frac{D}{M} &  \frac{\bone^\top}{M} & 0 \\
         (\diag(\btau))^{-1}\br & -(\diag(\btau))^{-1}  & (\diag(\btau))^{-1}\bk \\
	\frac{\beta}{\tau_z} & \frac{\bone^\top}{\tau_z} & -\frac{1}{\tau_z}\end{bmatrix}\begin{bmatrix}  {y}   \\
         {\bx}  \\
	{z}\end{bmatrix} \nonumber\\&\quad\quad\quad\quad+ \begin{bmatrix}  \bzero^\top &  \frac{1}{M}  \\
         (\diag(\btau))^{-1}(\bI-\bk\bone^{\top})  & 0 \\
	\bzero^\top & 0 \end{bmatrix}\begin{bmatrix}  {\bx^{\star}}  \\
	{d}\end{bmatrix}.
\end{align}
\begin{figure*}[b]
    \begin{equation} \label{eq:ssvalues}
\begin{aligned}
 \begin{bmatrix}
	y_{\mathrm{ss}}   \\
	\bx_{\mathrm{ss}} \\
	z_{\mathrm{ss}}
	\end{bmatrix} &=      \begin{bmatrix}
		\bone^\top\frac{\bDelta}{D}(\bI-\bk\bone^\top) &  -\bone^\top\frac{\bDelta}{D^2}\left(\br+\beta\bk\right) -\frac{1}{D} \\
		\bDelta(\bI-\bk\bone^\top) & -\frac{\bDelta}{D}\left(\br+\beta\bk\right) \\
	\left(\frac{\beta}{D}+1\right)\bone^\top\bDelta(\bI-\bk\bone^\top)  &  -\left(\frac{\beta}{D}+1\right)\bone^\top\frac{\bDelta}{D}\left(\br+\beta\bk\right)-\frac{\beta}{D}
	\end{bmatrix}\begin{bmatrix}
	\bx^{\star} \\
        d
	\end{bmatrix}\\
\text{where,~}	\bDelta &= \left(\bI-\left(\frac{\beta}{D}+1\right)\bk\bone^\top -\frac{\br}{D}\bone^\top\right)^{-1}.
\end{aligned}
    \end{equation}
\end{figure*}
Let us denote by ($\bx_{\mathrm{ss}}, y_{\mathrm{ss}},z_{\mathrm{ss}}$) the finite steady-state values of the respective quantities. Setting the time derivatives in~\eqref{eq:statespace} to zero and solving for the steady-state quantities, one obtains~\eqref{eq:ssvalues}. In certain instances, the time-constants associated with dynamics of $(\bx,y)$ are markedly smaller than those associated with the supervisory control dynamics of $z$; the distinction in power systems motivates separate analysis of primary- and secondary-control design. To delve deeper into such analysis, we quantify the equilibrium without secondary control next.

\subsection{Equilibrium without Secondary Control}\label{subsec:slow}
Consider a supply-demand system without secondary control. Specifically, deviating from~\eqref{eq:y_dynamics}-\eqref{eq:agc}, consider
\begin{subequations}\label{eq:primarydynamics}
\begin{align}
    M \dot{y} +  D y &= \bone^\top\bx - d,\\
    \diag(\btau)\dot{\bx} &= \bx^\star-\bx + \br y,
\end{align}
\end{subequations}
where the reference supply is the constant vector $\bx^\star$ from an optimal dispatch problem. In power system parlance, the above model represents a collection of coherent generators equipped with primary control, but no secondary control. It can also represent system behavior during the time-horizon over which secondary control action is slow to act. (See Fig.~\ref{fig:Block diag}.) The steady state of~\eqref{eq:primarydynamics} in such a setting is given by
\begin{subequations}\label{eq:primary}
\begin{align}
     y_\mrss&=\frac{\bone^\top\bx^\star-d}{D-\bone^\top\br},\label{seq:primary-y}\\
    \bx_\mrss&=\bx^\star+\frac{\br (\bone^\top\bx^\star-d)}{D-\bone^\top\br}.\label{seq:primary-x}
\end{align}
\end{subequations}

Interpreting~\eqref{eq:primary} as the steady state of a system without supervisory control tempts to answer the pedagogical power-system question: \emph{Is primary control adequate to ensure frequency restoration and serving unaccounted demands?} Specifically, is it possible to get $y_\mrss=0$ and/or $\bone^\top\bx_\mrss=d$. Reaffirming the conventional wisdom, \eqref{seq:primary-y} dictates that unless $\bone^\top\bx^\star=d$, (implying that the demand is precisely accounted for in the centralized dispatch~\eqref{eq:RA}) the frequency deviation from nominal is non-zero. For the second question, premultiplying~\eqref{seq:primary-x} by $\bone^\top$ and rearranging yields $$\bone^\top \bx_\mrss-d=(\bone^\top\bx^\star-d)\frac{D}{D-\bone^\top\br}.$$
Since $D\neq0$, and $\bone^\top\bx^\star\neq d$ when the complete demand is not accounted for, one sees that $\bone^\top\bx_\mrss\neq d$ for finite $\br$. 
 
\section{Design}\label{sec:design}
We have so far presented the system-, agent-, and supervisory-control dynamics for a supply-demand setup representative of power systems. To recapitulate, the workflow begins with determining an optimal supply dispatch $\bx^\star$ adequate to serve an anticipated demand $\hd$. For the demand realization to be $d$ in real time, the integrated system acquires a steady-state characterized by~\eqref{eq:lek1}; and features a quasi-steady state given by~\eqref{eq:primary}. Thus, in response to the change in demand from $\hd$ to $d$, the steady-state supply vector changes from $\bx^\star$ to $\bx^\star+\delta\bx$, where the supply adjustment $\delta\bx$ is precisely quantified by~\eqref{seq:lek1:xss} and~\eqref{seq:primary-x}. In principle, for the realized demand $d$, one could resolve~\eqref{eq:RA} to obtain a new cost optimal dispatch $\check{\bx}^\star$. However, for frequent demand changes such centralized interventions are often limited due to computational and communication restrictions. Thus, in the event of a demand change, the supply vector first adjusts per~\eqref{seq:primary-x}, then settles at~\eqref{seq:lek1:xss} (if equipped with supervisor control), before being updated to the new optimal dispatch $\check{\bx}^\star$.
\subsection{Secondary-control Participation Factors}
The supervisory control~\eqref{eq:agc} was motivated by the need to restore the imbalance metric such that $y_\mrss=0$, equivalently satisfying the supply-demand balance, i.e., $\bone^\top\bx_\mrss=d$. However, these attributes are not yet evident from the steady-state quantities computed in~\eqref{eq:ssvalues}. To elucidate that the supervisory control indeed yields $y_\mrss=0$, one needs to incorporate the prevalent design practice of enforcing $\bone^\top\bk=1$. In doing so, additional interesting aspects regarding the steady state of the integrated system show up as provided by the next result.
\begin{lemma}\label{le:k1}
Considering a well-defined parameterization of LTI~\eqref{eq:statespace}, if $\bone^\top\bk=1$, then the steady-state quantities are:
\begin{subequations}\label{eq:lek1}
\begin{align}
    y_\mrss&=0,\label{seq:lek1:yss}\\
    \bx_\mrss&=\bx^\star+\bk(d-\bone^\top\bx^\star),\label{seq:lek1:xss}\\
    z_\mrss&=d.\label{seq:lek1:zss}
\end{align}
\end{subequations}
\end{lemma}
\begin{IEEEproof}
The proof proceeds with re-evaluating the steady state identified in~\eqref{eq:ssvalues} while substituting $\bone^\top\bk=0$. To simplify notation, consider the definition
$$\bnu=-\frac{\br}{D}-\left(\frac{\beta}{D}+1\right)\bk,$$
that allows us to express
\begin{equation}\label{eq:bDelta}
    \bDelta=(\bI+\bnu\bone^\top)^{-1}=\bI-\frac{\bnu\bone^\top}{1+\bone^\top\bnu},
\end{equation}
where the second equality follows from the Woodbury matrix identity. Using $\bone^\top\bk=1$ it follows that
$$\bone^\top\bnu=-\frac{\bone^\top\br-\beta}{D}-1.$$
Substituting the above in~\eqref{eq:bDelta} reads
\begin{equation}
    \bDelta=\bI+\frac{D\bnu\bone^\top}{\bone^\top\br+\beta}.
\end{equation}
Replacing the above expression for $\bDelta$ in the coefficients for evaluation of $\bx_\mrss$ in~\eqref{eq:ssvalues}, one finds after some algebraic manipulations that
\begin{subequations}\label{eq:values}
\begin{align}
    \bDelta(\bI-\bk\bone^\top)&=(\bI-\bk\bone^\top)\\
    \frac{\bDelta}{D}\left(\br+\beta\bk\right)&=-\bk,
\end{align}
\end{subequations}
resulting in the steady-state supply vector being
\begin{equation}\label{eq:xss1}
    \bx_\mrss=\bx^\star+\bk(d-\bone^\top\bx^\star).
\end{equation}
From~\eqref{eq:xss1}, we get $\bone^\top\bx_\mrss=d$, thus providing~$y_\mrss=0$. The same can be alternatively verified by replacing~\eqref{eq:values} in~\eqref{eq:ssvalues}. From the above replacement, one further obtains $z_\mrss=d$
\end{IEEEproof}

\begin{remark}\longthmtitle{Interpreting the steady state characterized by Lemma~\ref{le:k1}}
In addition to showing that the considered integrated system dynamics of~\eqref{eq:statespace} can successfully restore $y_\mrss=0$, Lemma~\ref{le:k1} provides interesting insights on $\bx_\mrss$ and $z_\mrss$. Specifically, while demand $d$ is an exogenous variable which may not be directly measurable, the supervisory control state $z$ tracks $d$ in steady state. Moreover, during steady state, the excess demand to be served $d-\bone^\top\bx^\star$ (alternatively interpreted as the unaccounted demand in solving optimal dispatch~\eqref{eq:RA}) is allocated to different agents in the proportion determined by $\bk$. Since the entries of $\bk$ sum up to unity and determine the extent to which each supplier participates in picking up the excess demand, these are oftentimes referred to as \emph{participation factors}. \oprocend
\end{remark}

While Lemma~\ref{le:k1} reaffirms the need to have $\bone^\top \bk = 1$, a pertinent follow-up consideration relates to the choice of individual entries of $\bk$. It turns out that one particular option ensures the adjusted steady-state supply per~\eqref{seq:lek1:xss} coincides with the optimal allocation as per~\eqref{eq:RA}, i.e.,
\begin{equation}\label{eq:q1}
\bx^\star+\bk(d-\bone^\top\bx^\star)=\check{\bx}^\star.
\end{equation}
Consider the following reasonable assumptions:
\begin{assumption}\label{as:nonbinding}
 None of the agents are at their limits, before and after adjustments catering demand variation. Specifically, inequalities in~\eqref{eq:limitx} are non-binding.
\end{assumption}
\begin{assumption}\label{as:quadratic}
Cost functions $f_n(\cdot)$ are strictly convex quadratic functions. Stated alternatively,
$$f_n(\chi_n)=a_n\chi_n^2+b_n\chi_n+c_n,~\forall n,$$
where $a_n>0$.
\end{assumption}
With the aid of Karuhn-Kush-Tucker (KKT) conditions that are satisfied by minimizers of~\eqref{eq:RA}, and under assumptions~\ref{as:nonbinding} and~\ref{as:quadratic}, optimal sharing to ensure~\eqref{eq:q1} can be obtained with
\begin{equation}\label{eq:optimal_rho}
    \bk=\frac{\bA^{-1}\bone}{\bone^\top\bA^{-1}\bone},
\end{equation}
where $\bA:=\diag(\{a_n\}_{n=1}^N)$~\cite{Baros21AGC}. Note that since \eqref{eq:RA} is a strictly convex problem under assumption~\ref{as:quadratic}, it features a unique minimizer for a given demand $d$. Thus, the participation vector $\bk$ identified in~\eqref{eq:optimal_rho} is a unique solution to the desired optimality requirement.

\subsection{Primary-control Regulation Gains}
In systems not equipped with secondary control, the steady-state supply depends on the primary-control regulation gains~$\br$, as characterized in~\eqref{eq:primary}. In such cases, it is often sought that the change in the steady-state supply in~\eqref{seq:primary-x} by all agents is in proportion to their maximum capacity. Specifically, for what choice of $\br$ does the following hold:
\begin{equation}\label{eq:q2}
\frac{\delta x_i}{\delta x_j}=\frac{\overline{x}_i}{\overline{x}_j},~\forall i,j\in\{1,\dots,N\},\end{equation} where $$\delta\bx=\frac{\br (\bone^\top\bx^\star-d)}{D-\bone^\top\br}.$$

The aforementioned pursuit tacitly incorporates assumption~\ref{as:nonbinding}. Since $\delta\bx$ is defined to be proportional to $\br$, the proportional sharing requirement in~\eqref{eq:q2} requires 
\begin{equation}\label{eq:proportional_r}
    \br=\kappa\overline{\bx},
\end{equation}
where $\kappa$ is a proportionality constant. While~\eqref{eq:proportional_r} suffices to ensure proportionate sharing, it allows an ambiguity in the choice of $\kappa$, indicating a non-uniqueness in the choice of $\br$. Indeed, the choice of $\kappa$ determines the deviation of $y_\mrss$ from zero; cf.~\eqref{seq:primary-y}. The choice of $\br$ can be uniquely made by imposing an additional requirement on the desired \emph{level of regulation}, which is specified as
\begin{equation}\label{eq:regulation}
    \frac{\bone^\top\bx^\star-d}{y_\mrss}=\pi.
\end{equation}
For a given $\pi$, the value of $\kappa$ can be obtained by replacing~\eqref{eq:proportional_r} in~\eqref{seq:primary-y} as
\begin{equation}
    \kappa=\frac{D-\pi}{\bone^\top\overline{x}}.
\end{equation}
Notably, this matches the result in~\cite{Guggilam-2018}. 

\section{Conclusion}\label{sec:conclusion}
This work puts forth an abstract power-system representation that is rich enough to succinctly capture primary-, secondary-, and tertiary-control schemes central to system operation. The developed model is then used to characterize system steady states, under different configurations, that conspicuously portray the role of key design parameters. Such characterization allows justifying conventional design choices encountered in classical power-system literature. It is envisaged that the developed framework will enable analysis and design of advanced control for future power grids. 
\bibliographystyle{IEEEtran}
\bibliography{Ref}
\end{document}